\documentclass{amsproc}
\usepackage{graphicx}
\usepackage{amssymb}
\usepackage{epstopdf}
\usepackage{amsmath,amsthm,amscd,amssymb}
\usepackage{latexsym}
\usepackage[colorlinks,citecolor=red,pagebackref,hypertexnames=false]{hyperref}
\usepackage{geometry}                
\numberwithin{equation}{section}

\theoremstyle{plain}
\newtheorem{theorem}{Theorem}[section]
\newtheorem{lemma}[theorem]{Lemma}
\newtheorem{corollary}[theorem]{Corollary}

\theoremstyle{definition}
\newtheorem{definition}[theorem]{Definition}

\newtheorem{case[theorem]}{Case}

\newcommand{\ag}{}

\def\dH{\dim_{{\mathcal H}}}
\def\R{\Bbb R}

\theoremstyle{remark}

\numberwithin{equation}{section}

\begin{document}

\title[Three-point configurations determined by subsets of the Euclidean plane]{\parbox{14cm}{\centering{Three-point configurations determined by subsets of the  plane, \\
a bilinear operator and applications  to discrete geometry}}}


\author{Allan Greenleaf and Alex Iosevich}
\address{Department of Mathematics \\ University of Rochester\\ Rochester, NY 14627}
\email{allan@math.rochester.edu, iosevich@math.rochester.edu}

\thanks{The authors were supported by NSF grants DMS-0853892 and DMS-1045404.}

\begin{abstract} We prove that if the Hausdorff dimension of a compact set $E \subset {\Bbb R}^2$ is greater than \nolinebreak$\frac{7}{4}$, then
the set of {\ag three-point configurations determined by $E$ has positive three-dimensional  measure}. We 
establish this  by showing that {\ag a} natural measure on the set of {\ag such configurations} 
has {\ag Radon-Nikodym derivative} in $L^{\infty}$  if $\dH(E)>\frac{7}{4}$, 
and the index $\frac{7}{4}$ in this last result cannot, in general, be improved. 
This problem naturally leads to the study of a bilinear convolution operator,
$$ B(f,g)(x)=\int \int f(x-u) g(x-v)\, dK(u,v),$$ where $K$ is surface measure on the set $ \{(u, v) \in\R^2 \times \R^2: |u|=|v|=|u-v|=1\}$, and  we prove a scale of estimates  that includes
$B:L^2_{-\frac{1}{2}}({\Bbb R}^2) \times L^2({\Bbb R}^2) \to L^1({\Bbb R}^2)$ on positive functions. 

As an application of our main result, it follows that {\ag for finite sets  of cardinality $n$   and belonging to a natural class of discrete sets in the plane}, the maximum number of times a given three-point configuration arises is $O(n^{\frac{9}{7}+\epsilon})$ (up to congruence), improving upon the known bound of $O(n^{\frac{4}{3}})$ in this context. 
\end{abstract} 

\maketitle


\section{Introduction}

\vskip.125in 

The classical Falconer distance conjecture says that if 
a compact set $E \subset {\Bbb R}^d$, $d \ge 2$,
has Hausdorff dimension $dim_{{\mathcal H}}(E)>\frac{d}2$,
then the one-dimensional Lebesgue measure ${\mathcal L}^1(\Delta(E))$ of its {\it distance set},
$$ \Delta(E):=\{|x-y|\in\Bbb R: x,y \in E\},$$ 
is positive. Here, and throughout, $|\cdot|$ denotes the Euclidean distance. A beautiful example due to Falconer, based on the integer lattice, shows that the exponent $\frac{d}{2}$ is best possible. The best results currently known, culminating almost three decades of efforts by Falconer \cite{Fal86}, Mattila \cite{Mat87}, Bourgain \cite{B94} and others, are due to Wolff \cite{W99} for $d=2$  and  Erdo\v{g}an \cite{Erd05} for $d\ge 3$. They prove that 
${\mathcal L}^1(\Delta(E))>0$ if 
$$dim_{{\mathcal H}}(E)>\frac{d}{2}+\frac{1}{3}.$$ 

Since  two-point configurations are equivalent, up to Euclidian motions of $\R^d$, precisely if the corresponding distances are  the same, one may think of the Falconer conjecture as stating that the set of  two-point configurations determined by a compact $E$ of sufficiently high Hausdorff dimension has positive measure. 
A natural extension of the Falconer problem is then the question: 
\medskip

\centerline{{\bf Q:} \it For   $N\ge 3$, how great does the Hausdorff dimension of 
a compact set need to be, }
\centerline{ \it in order to ensure that the  set of $N$-point configurations it determines is of positive measure?}
\medskip

To make this more precise, define 
the {\it space of $(k+1)$-point configurations in $E$}, or
the {\it quotient space of (possibly degenerate)  $k$-simplices with vertices in $E$}, modulo Euclidian motions, as
$$ T_k(E): =E^{k+1}/ \sim,$$ where $E^{k+1}=E \times E \times \dots \times E$ ($k+1$ times) and the congruence relation
$$ (x^1, x^2, \dots, x^{k+1}) \sim (y^1, y^2, \dots, y^{k+1})$$
holds  iff there exists an element $R$ of the orthogonal group $O(d)$ and a translation $\tau \in {\Bbb R}^d$ such that 
$$ y^j=\tau+R(x^j),\quad 1\le j\le k+1.$$  
Observe that we may identify $T_k(E)$ as a subset of ${\Bbb R}^{{k+1 \choose 2}}$, since rigid motions may be encoded by fixing distances, and this induces ${k+1 \choose 2}$-dimensional Lebesgue measure on $T_k(E)$. 
The problem under consideration was first taken up in \cite{EHI11} where it was shown that 
$$\hbox{ if }dim_{{\mathcal H}}(E)>\frac{d+k+1}{2}, \ \text{then} \  {\mathcal L}^{{k+1 \choose 2}}(T_k(E))>0.$$ 
Unfortunately, these results do not give a non-trivial exponent for what are arguably the most natural cases, namely {\ag three-point configurations in $\R^2$, four-point configurations in $\R^3$ and, more generally,  $(d+1)$-point configurations (generically spanning $d$-dimensional simplices)} in ${\Bbb R}^d$. 
(Nor does it yield results for $(d-1)$-simplices.)
Here, we partially fill this gap by establishing a non-trivial exponent for three-point configurations in the plane. 
\medskip

As for counterexamples, it is easy to see that ${\mathcal L}^{{k+1 \choose 2}}(T_k(E))>0$ does not hold if the Hausdorff dimension of $E$ is less than or equal to $d-1$; to see this, just take $E$ to be a subset of a $(d-1)$-dimensional plane. We do not currently know if more restrictive {\ag conditions} exist in this context. However, more restrictive counterexamples do exist if we consider the following related question. For $t=\{t_{ij} \}$ any symmetric matrix with zeros on the diagonal, let

$$ {\mathcal S}_t^k(E)=\big\{(x^1, \dots, x^{k+1}) \in E^{k+1}: |x^i-x^j|=t_{ij},\,\forall i,j \, \big\}.$$ 
Conditions under which 

\begin{equation} \label{sard} \dim_{{\mathcal M}}\left({\mathcal S}_t^k\left(E\right)\right) \leq (k+1) \dim_{{\mathcal H}}(E)-{k+1 \choose 2}=(k+1)\left(\dH\left(E\right)-\frac{k}2\right),
\end{equation} where $\dim_{{\mathcal M}}$ denotes the Minkowski dimension, are analyzed in \cite{EIT11} in the case $k=1$ in a rather general setting and in \cite{CEHIT10} in the case $k>1$. (See \cite{Falc86,M95} for background on $\dH, \dim_{{\mathcal M}}$ and connections with harmonic analysis.) 

The estimate (\ref{sard}) follows easily if one can show that 

\begin{equation} \label{bigmama} (\mu \times \mu \times \dots \times \mu) \big\{(x^1, \dots, x^{k+1}): t_{ij} \leq |x^i-x^j| \leq t_{ij}+\epsilon,\,\forall i,j\, \big\} \lesssim \epsilon^{k+1 \choose 2}\hbox{ as }\epsilon\searrow0,\end{equation} 
where $\mu$ is a Frostman measure (defined in (\ref{frostman}) below) on $E$, under the assumption that 
$\dH(E)>s_0$ for some
threshold  $s_0<d$. 
This is shown in \cite{CEHIT10} under the assumption that the Hausdorff dimension of $E$ is greater than $\frac{k}{k+1}d+\frac{k}{2}$, but  observe that this only yields a non-trivial exponent (less than $d$) if ${k \choose 2}<d$ and, in particular, does not cover the important case of $k=d$. 

Our main 
result is the following. 
\begin{theorem} \label{main}  Let $E \subset {[0,1]}^2$ be compact and $\mu$ a Frostman measure on $E$.
\medskip

(i) If   $\dH(E)>\frac{7}{4}$,  then  estimate (\ref{bigmama}) holds with  $d=k=2$. 

\vskip.125in 

(ii)  If $\dH(E)>\frac{7}{4}$ , then  $ {\mathcal L}^3(T_2(E))>0$. 
\end{theorem} 

\vskip.125in

The proof that  part (i) of Theorem \ref{main} implies part (ii) is presented in Sec. \ref{reduction} below;  part (i) is then proved by analysis of a bilinear operator (or trilinear form) in Sections \ref{reduction}, \ref{tritobi},  and \ref{bi}. 

We observe that the result in part (i)  is sharp in the following sense. Define a measure $\nu$ on $T_2(E)$ by the relation 

\begin{equation} \label{nu} \int\! f(t_{12}, t_{13}, t_{23}) d\nu(t_{12}, t_{13}, t_{23})=\int\!\int\! \int\! f \left(|x^1-x^2|, |x^1-x^3|, |x^2-x^3| \right) d\mu(x^1) d\mu(x^2) d\mu(x^3),\end{equation} 
where $\mu$ is any Frostman measure on $E$. We shall prove that the Radon-Nikodym derivative $\frac{d\nu}{dt} \in L^{\infty}$, which is just a rephrasing of the statement that (\ref{bigmama}) holds for $d=k=2$,  if the Hausdorff dimension of $E$ is greater than $\frac{7}{4}$. 
On the other hand, we also use a variant of Mattila's example from \cite{Mat87} to show that if $s<\frac{7}{4}$, then $\frac{d\nu}{dt}$ need not be, in general, in $L^{\infty}$, in the sense that for every $s<\frac{7}{4}$ there exists a set $E$ of Hausdorff dimension $s$ and a Frostman measure $\mu$ supported on $E$, such that $\frac{d\nu}{dt}\not\in L^{\infty}$. (This issue is taken up in  Sec. \ref{sharpness}.) Thus, in order to try to improve part (ii) of the Theorem, i.e., to prove that ${\mathcal L}^3(T_2(E))>0$ if  $\dH(E)=s_0$, for some $s_0\le \frac{7}{4}$, it would be reasonable to try to obtain an $L^p$, rather an $L^{\infty}$ bound on the measure $\nu$ defined by (\ref{nu}). We hope to address this  in a subsequent paper. 

\medskip

Theorem \ref{main} may be viewed as a local version of the following theorem due to Furstenberg, Katznelson and Weiss \cite{FKW90}; see also \cite{B86,Z06} for subsequent results along these lines. 
\begin{theorem} \label{fkw90II} \cite{FKW90} Let $E \subset {\Bbb R}^2$ be of positive upper Lebesgue density, in the sense that 

$$ \limsup_{R \to \infty} \frac{{\mathcal L}^d \{E \cap {[-R,R]}^2 \}}{{(2R)}^2}>0,$$ 
where ${\mathcal L}^2$ denotes $2$-dimensional Lebesgue measure. For $\delta>0$, let $E_{\delta}$ denote the $\delta$-neighborhood of $E$. Then, given vectors $u,v$ in ${\Bbb R}^2$, there exists $l_0$ such that for any  $l>l_0$ and  $\delta>0$, there exist $x,y,z\in E_{\delta}$ forming a triangle congruent to $\{{\bf 0},lu,lv\}$, where ${\bf 0}$ denotes the origin in ${\Bbb R}^2$. 
\end{theorem} 

We note in passing that it is generally believed that the conclusion of Theorem \ref{fkw90II} still holds if the $\delta$-neighborhood of $E$ is replaced by $E$ under an additional assumption that the triangles under consideration are non-degenerate. 
For degenerate triangles, i.e., allowing line segments, the necessity of considering the $\delta$-neighborhood of $E$ was established by Bourgain (see \cite{FKW90}). 

In contrast to Theorem \ref{fkw90II}, we are able in the local version to go beyond subsets of the plane of positive Lebesgue measure, and we do not need to allow for dilations of the triangles. On the other hand, we only obtain a positive Lebesgue measure's worth of the possible three-point configurations,  not all of them.

It is also not difficult to show (see Sec. \ref{reduction}) that if the estimate (\ref{bigmama}) holds under the assumption that 
$\dH(E)>s_0$,
then ${\mathcal L}^{k+1 \choose 2}(T_k(E))>0$ for 
these sets.
In \cite{CEHIT10}, a number of estimates of the type (\ref{bigmama}) are proved, but, as we note above,  do not cover the cases $k=d$ or $k=d-1$.

\vskip.125in 

\subsection{A combinatorial perspective} Finite configuration problems have their roots in geometric combinatorics. For example, the Falconer distance problem is a continuous analog of the celebrated Erd\H os distance problem; see \cite{SolToth01,KT04,BMP05,Sz97} and the references  therein. The discrete precursor of the problem discussed in this paper is the following question posed by Erd\H os and Purdy (see \cite{BMP05,EP71},  and also \cite{EP75,EP76,EP77,EP78,EP95}): 

\vskip.125in 

{\bf Q:} \emph{What is the maximum number of mutually congruent $k$-simplices with vertices from among a set of $n$ points in ${\Bbb R}^d$? }

\vskip.125in

In Sec. \ref{application} we shall see that Theorem \ref{main} (ii) implies that for a large class of finite sets  $P$ of cardinality $n$ in $\R^2$,   namely those that are {\it $s$-adaptable}, the maximum number of mutually congruent triangles determined by points of $P$ is $O(n^{\frac{9}{7}+\epsilon})$. 

For explicit quantitative connections between discrete and continuous finite configuration problems in other contexts, see, for example, \cite{HI05}, \cite{IL05} and \cite{IRU10}. 

\subsection{Notation:} Throughout the paper, $X \lesssim Y$ means that there exists $C>0$ such that $X \leq CY$ and $X \approx Y$ means that $X \lesssim Y$ and $Y \lesssim X$. We also define $X \lessapprox Y$ as follows. If $X$ and $Y$ are quantities that depend on a large parameter $N$, then $X \lessapprox Y$ means that for every $\epsilon>0$ there exists $C_{\epsilon}>0$ such that $X \leq C_{\epsilon}N^{\epsilon}Y$, while if $X$ and $Y$ depend on a small parameter $\delta$, then $X \lessapprox Y$ means that for every $\epsilon>0$ there exists $C_{\epsilon}>0$ such that $X \leq C_{\epsilon} \delta^{-\epsilon}Y$ as $\delta$ tends to $0$. 
\medskip

\section{Reduction of the proof  to the estimation of a tri-linear form} 
\label{reduction} 

\vskip.125in 

We shall work exclusively with  Frostman measures. Recall that a probability measure $\mu$ on a compact set $E\subset\R^d$ is a {\it Frostman measure}  if, for any ball $B_{\delta}$ of radius $\delta$, 
\begin{equation} \label{frostman} \mu(B_{\delta}) \lessapprox \delta^s, \end{equation} 
where $s=\dH(E)$. For discussion and  proof of the existence of such measures see, e.g., \cite{M95}. 
\medskip

Let $\mu$ be a Frostman measure on $E$. Cover $T_2(E)$ by cubes of the form 
$$(t^l_{12}-\epsilon_l, t^l_{12}+\epsilon_l) \times  (t^l_{13}-\epsilon_l, t^l_{13}+\epsilon_l) \times  (t^l_{23}-\epsilon_l, t^l_{23}+\epsilon_l).$$ It follows that 

\begin{equation} \label{pigeononparade} 1=(\mu \times \mu \times \mu) \big\{E \times E \times E\big\} \leq \sum_l (\mu \times \mu \times \mu) \big\{(x^1,x^2,x^3): t^l_{ij}-\epsilon_l \leq |x^i-x^j| \leq 
t^l_{ij}+\epsilon_l,\, \forall i,j\,\big \}. \end{equation}

Suppose that we could show that this expression is  $\lesssim \sum \epsilon_l^3$.
It would then follow, by definition of sets of measure $0$, that the three dimensional Lebesgue measure of $T_2(E)$ is positive. 

In light of (\ref{frostman}), to establish the positive measure of $T_2(E)$ we may assume that $t_{ij} \ge c>0$. To see this, observe that if each $t_{ij}$ is $\leq r$, then fixing $x^1$ results in $x^2$ and $x^3$ being contained in a ball of radius $r$ centered at $x^1$. It follows that 
$$ (\mu \times \mu \times \mu) \big\{E \times E \times E\big\} \leq \sum_l (\mu \times \mu \times \mu) \big\{(x^1,x^2,x^3): t^l_{ij}-\epsilon_l \leq |x^i-x^j| \leq  t^l_{ij}+\epsilon_l,\, \forall i,j\,\big \} \leq Cr^{2s},$$ and taking $r$ to be small enough, this expression is $\leq \frac{1}{10}$. This means that in place of equality on the left hand side of (\ref{pigeononparade}), we have an inequality with $1$ replaced by $\frac{9}{10}$ and the rest of the argument goes through as before. 

Therefore the proof of Theorem \ref{main} (i) is reduced to proving the trilinear estimate 

\begin{equation} \label{mamaexpression} 
\Lambda^{\epsilon}_t(\mu, \mu, \mu) := \int \int \int \sigma_{t_{12}}^{\epsilon}(x^1-x^2) \sigma_{t_{13}}^{\epsilon}(x^1-x^3) \sigma_{t_{23}}^{\epsilon}(x^2-x^3) d\mu(x^1) d\mu(x^2) d\mu(x^3) \lesssim 1.\end{equation}
Here, $t=(t_{12}, t_{13}, t_{23})$, $\sigma_r$ is arc length measure on the circle of radius $r$ in $\R^2$, and $\sigma_r^\epsilon=\sigma_r*\rho_\epsilon$, where $\rho_\epsilon(x)=\epsilon^{-2}\rho(\frac{x}\epsilon)$ is an approximate identity with $\rho\in C_0^\infty(\{|x|\le 1\})$, $\rho\ge 0$, $\int \rho(x) dx=1$. Note that the right hand side is $1$ instead of $\epsilon^3$ because the characteristic function of the annulus of radius $t_{ij}$ and thickness $\epsilon$, divided by $\epsilon$, is dominated by $\sigma_{t_{ij}}^{\epsilon}$. We now turn to the proof of (\ref{mamaexpression}).

\vskip.125in 

\section{Reducing the  trilinear form estimate   to a bilinear operator estimate} 
\label{tritobi} 

\vskip.125in

Define trilinear forms
\begin{equation} \label{analyticfamily} 
\Lambda^{\epsilon}_t(f_1, f_2, f_3):= \int \int \int \sigma_{t_{12}}^{\epsilon}(x^1-x^2) \sigma_{t_{13}}^{\epsilon}(x^1-x^3) \sigma_{t_{23}}^{\epsilon}(x^2-x^3) f_1(x^1) f_2(x^2) f_3(x^3)dx^1\, dx^2\, dx^3, \end{equation} 
and consider $\Lambda_t^\epsilon(\mu_{-\alpha}^\delta, \mu,\mu_{\alpha}^\delta)$, where 
\begin{equation} \label{spongebob} 
\mu_{\alpha}(x):=\frac{2^{\frac{2-\alpha}{2}}}{\Gamma(\alpha/2)}\big( \mu* {| \cdot |}^{-2+\alpha}\big)(x),
\end{equation} 
initially defined for $Re(\alpha)>0$, is extended to the complex plane by analytic continuation, and
$$ \mu^{\delta}(x):=\mu* \rho_{\delta}(x),$$ 
and $\rho_{\delta}(x)=\delta^{-2} \rho(x/ \delta)$ is an approximate identity as above. Observe that 
$ \widehat{\mu}_{\alpha}^\delta(\xi)=C_{\alpha} \widehat{\mu}(\xi) \widehat{\rho}(\delta \xi)  {|\xi|}^{-\alpha}$, where 
\begin{equation}\label{blowup}
C_{\alpha} = \frac{2\pi \cdot 2^{\frac{\alpha}{2}}}{\Gamma(\frac{2-\alpha}{2})}.
\end{equation}
(See \cite{GS58} for relevant calculations.)
This shows, in view of Plancherel, that $\mu_{\alpha}^\delta$ is an $L^2({\Bbb R}^2)$ function, with bounds depending on 
$\delta$. Moreover, since we have compact support, this shows that one has a trivial finite upper bound on the trilinear form with constants depending on $\delta$. Taking the modulus in (\ref{spongebob}), we see that
$$|\mu_\alpha^\delta(x)|\leq \left| \frac{2^{\frac{2-\alpha}{2}}}{\Gamma\left(\alpha/2 \right)} \right| (\mu^{\delta} * |\cdot |^{-2+\text{Re}(\alpha)})(x) = 2^{\frac{2-\text{Re}(\alpha)}{2}}\frac{\Gamma(\text{Re}(\alpha)/2)}{|\Gamma(\alpha/2)|}\mu_{\text{Re}(\alpha)}^\delta(x)$$
and note that the right hand side is non-negative.

Now define 
\begin{equation} \label{threelinesexpression} F(\alpha):=\Lambda^{\epsilon}_t(\mu_{-\alpha}^\delta, \mu, \mu_{\alpha}^\delta)= \langle B(\mu^{\delta}_{-\alpha}, \mu^{\delta}), \mu^{\delta}_{\alpha}\rangle,\end{equation} 
where $\langle \cdot, \cdot\rangle$ is the $L^2({\Bbb R}^2)$ inner product and $B$ is the bilinear operator given by the relation 
\begin{equation} \label{mainoperator} B^{\epsilon}(f,g)(x):=\int \int f(x-u) g(x-v) \sigma_a^{\epsilon}(u) \sigma_b^{\epsilon}(v) \sigma^{\epsilon}(u-v)\, du\, dv. \end{equation} 
Here, for simplicity we have rescaled one  side of the triangle to have unit length;
 the other two,  $a,b\lesssim 1$, are bounded away from 0.

Our main bilinear estimate is the following, which is proved in \S\ref{bi}.
\begin{theorem} \label{bilineartheorem} Let $B^{\epsilon}$ be defined as above and suppose that $f,g \ge 0$. Then 
\begin{equation} \label{bilinearestimate} B^\epsilon: L^2_{-\beta_1}({\Bbb R}^2) \times L^2_{-\beta_2}({\Bbb R}^2) \to L^1({\Bbb R}^2) \ \text{if} \ \beta_1+\beta_2=\frac{1}{2}, \beta_1, \beta_2 \ge 0 \end{equation} with constants independent of $\epsilon$. 
\end{theorem} 

Using (\ref{bilinearestimate}), we see that, with $F(\alpha)$ defined as in (\ref{threelinesexpression}), we have 
\begin{equation} \label{battlebegins} 
|F(\alpha)| \lesssim \  \langle B^\epsilon(\mu^{\delta}_{-Re(\alpha)}, \mu^{\delta}\rangle  \leq {||B^\epsilon(\mu^{\delta}_{-Re(\alpha)})||}_{L^1({\Bbb R}^2)} \cdot {||\mu^{\delta}_{Re(\alpha)}||}_{L^{\infty}({\Bbb R}^2)},
\end{equation} where $\lesssim$ symbol includes factors of the Gamma functions. 

\begin{lemma} \label{Linfinity} Suppose that $\mu$ is a Frostman measure on a set of Hausdorff dimension $>\frac{7}{4}$. Then 
$$ {||\mu_{\alpha}^{\delta}||}_{\infty} \lessapprox 1 \ \text{if} \ Re(\alpha)=\frac{1}{4}.$$ 
\end{lemma} 

To prove the lemma,  observe that if $Re(\alpha)=\frac{1}{4}$, 

\begin{equation}\nonumber
\mu_{\alpha}^{\delta}(x) \leq \int {|x-y|}^{-2+\frac{1}{4}} d\mu^{\delta}(y) \approx \sum_m 2^{m(2-\frac{1}{4})} \int_{|x-y| \approx 2^{-m}} d\mu^{\delta}(y) \lesssim  \sum_m 2^{m(2-\frac{1}{4})} 2^{-ms},
\end{equation} 
and this  is $\lessapprox 1$, since $\mu$ is a Frostman measure on a set of Hausdorff dimension $>\frac{7}{4}$. 
Substituting this into (\ref{battlebegins})
and applying 
(\ref{bilinearestimate}) with $\beta_1=\frac{3}{8}, \beta_2=\frac{1}{8}$, we see that,
if $Re(\alpha)=\frac{1}{4}$, 
\begin{equation}\label{product}
|F(\alpha)| \leq {||B^\epsilon(\mu^{\delta}_{-\frac{1}{4}}, \mu^{\delta})||}_{L^1({\Bbb R}^2)}
\lesssim {||\mu^{\delta}_{-\frac{1}{4}}||}_{L^2_{-\frac{3}{8}}({\Bbb R}^2)} \cdot {||\mu^{\delta}||}_{L^2_{-\frac{1}{8}}({\Bbb R}^2)}.
\end{equation}
A straightforward calculation using the definition of $\mu^{\delta}_{\alpha}$ from above shows that the square of either of the terms in (\ref{product}) is bounded by 
$$ \int \int {|x-y|}^{-\frac{7}{4}} d\mu(x) d\mu(y), $$ which is the energy integral of $\mu$ of order $\frac{7}{4}$. This integral is bounded since the Hausdorff dimension of $E$ is greater than $\frac{7}{4}$ and $\mu$ is a Frostman measure; see, e.g., \cite{Falc86,M95}. 

By symmetry, the same bound holds when $Re(\alpha)=-\frac{1}{4}$, because we can reverse the roles of $d\mu(x^1)$ and $d\mu(x^3)$. When $-\frac{1}{4}<Re(\alpha)<\frac{1}{4}$, we use the fact that $|F(\alpha)|$ is bounded from above with constants depending on $\delta$ as we noted in the beginning of this section. By the three lines lemma, (see, for example, I.I. Hirschman's version in \cite{H52}), we conclude that $\Lambda^{\epsilon}_t(\mu, \mu, \mu) \lesssim 1$, which completes the proof of Theorem \ref{main}, conditional on Theorem \ref{bilineartheorem}, which we now prove. 


\section{Estimating the bilinear operator } 
\label{bi}

Since we are assuming $f,g \ge 0$, we have 
\begin{equation} \label{L1initial} {||B^{\epsilon}(f,g)||}_{L^1({\Bbb R}^2)}=\int \int \int f(x-u)g(x-v) K^{\epsilon}(u,v)\, du\, dv\, dx, \end{equation} where 
$$ K^{\epsilon}(u,v)=\sigma_a^{\epsilon}(u)\sigma_b^{\epsilon}(v)\sigma^{\epsilon}(u-v);$$ 
recall that we scaled one of the sigmas to the unit radius. 
Let $\psi\in C_0^\infty(\{|x|\le 2\}),\, \psi\ge 0,\, \psi\equiv 1\hbox{ on }\{|x|\le 1\}$.
Then, it suffices to estimate 
$$ \int \int \int f(x-u)g(x-v) K^{\epsilon}(u,v)\, du\, dv\, \psi(x/R) dx$$ with bounds independent of $R \ge 1$. 
Using Fourier inversion, the expression (\ref{L1initial}) equals 
\begin{equation} \label{L1fourier} R^2 \int \int \widehat{f}(\xi) \widehat{g}(\eta) \widehat{K^{\epsilon}}(\xi, \eta) \widehat{\psi}(R(\xi+\eta)) d\xi d\eta. \end{equation}

We shall need the following stationary phase calculation. 

\begin{lemma} \label{stationaryphase} Let $K(u,v)=K^0(u,v)$, interpreted in the sense of distributions. We have 
\begin{equation} \label{twohands} \widehat{K}(\xi, \eta)=\widehat{\sigma}(U_{a,b}(\xi, \eta)), \end{equation} where 
$$ U_{a,b}: {\Bbb R}^4 \to {\Bbb R}^2$$ and is defined by 
\begin{equation} \label{mapping} U_{a,b}(\xi, \eta)=\left( a\xi_1+b \eta_1 \gamma_{a,b}+b \eta_2 \sqrt{1-\gamma^2_{a,b}}, 
a\xi_2-b \eta_1 \sqrt{1-\gamma^2_{a,b}}+b \gamma_{a,b} \eta_2 \right) \end{equation} with $ \gamma_{a,b}=\frac{a^2+b^2-1}{2ab}$. Consequently, 
\medskip

\begin{equation} \label{handsandfeet} \left|\widehat{K}^{\epsilon}(\xi, \eta)\widehat{\psi}(R(\xi+\eta))\right| \lesssim {(1+|\xi|+|\eta|)}^{-\frac{1}{2}}\end{equation} 
uniformly for $R \ge 1$. 
\end{lemma} 

\vskip.125in 

Recalling that, by the  method of stationary phase (see e.g. \cite{So93}, \cite{St93}), 
$$ |\widehat{\sigma}(\xi)| \lesssim {(1+|\xi|)}^{-\frac{1}{2}},$$ 
one sees that (\ref{handsandfeet})  will immediately follow from (\ref{twohands}),(\ref{mapping}).

To prove the lemma, parameterize the Cartesian product of two circles as 

$$ \{(a\cos(\theta), a\sin(\theta), b\cos(\phi),b\sin(\phi)) \}.$$ The restriction imposed by $\sigma(u-v)$ says that 

$$ dist((a\cos(\theta), a\sin(\theta)), (b\cos(\phi), b\sin(\phi))=1,$$ which implies via standard trigonometric identities that 

$$\cos(\theta-\phi)=\frac{a^2+b^2-1}{2ab} \equiv \gamma_{a,b}$$ 
and thus  $\theta-\phi=\pm \theta_{a,b}=\cos^{-1}(\gamma_{a,b})$. It follows that 

\begin{equation}\nonumber
\widehat{K}(\xi, \eta)=\int_0^{2 \pi} e^{2 \pi i (a\cos(\theta)\xi_1+a\sin(\theta)\xi_2+b\cos(\theta+\theta_{a,b})\eta_1+b\sin(\theta+\theta_{a,b})\eta_2)} d\theta=\widehat{\sigma}(U_{a,b}(\xi, \eta)),
\end{equation} 
as claimed. This proves (\ref{twohands}). The estimate (\ref{handsandfeet}) follows in the same way since $\sigma^{\epsilon}_a(x)=\sigma_a*\rho_{\epsilon}(x)$. 

Using Lemma \ref{stationaryphase}, Cauchy-Schwarz, and the assumption $\beta_1+\beta_2=\frac{1}{2}$, $\beta_1, \beta_2 \ge 0$, we estimate the square of (\ref{L1fourier}) by 
$$ \int {|\widehat{f}(\xi)|}^2 \left\{ R^2 \int {|\widehat{K^{\epsilon}}(\xi, \eta)|}^{4\beta_1} |\widehat{\psi}(R(\xi+\eta))| d\eta \right\} d\xi \cdot \int {|\widehat{g}(\eta)|}^2 \left\{ R^2 \int {|\widehat{K^{\epsilon}}(\xi, \eta)|}^{4\beta_2} |\widehat{\psi}(R(\xi+\eta))| d\xi \right\} d\eta$$
$$ \lesssim \int {|\widehat{f}(\xi)|}^2 {(1+|\xi|)}^{-2\beta_1} d\xi \cdot \int {|\widehat{g}(\eta)|}^2 {(1+|\eta|)}^{-2\beta_2} d\eta$$
$$={||f||}_{L^2_{-\beta_1}({\Bbb R}^2)}^2 \cdot {||g||}_{L^2_{-\beta_2}({\Bbb R}^2)}^2,$$ as desired, completing the proof of Theorem \ref{bilineartheorem} and thus the proof of Theorem \ref{main}. 

\vskip.125in 

\section{Sharpness of the trilinear estimate (\ref{mamaexpression})} 
\label{sharpness} 

\vskip.125in 

To understand the extent to which this result is sharp, we use a variant of the construction due to Mattila obtained for the case $k=1, d=2$ in \cite{Mat87}. See \cite{IS10,CEHIT10} where this issue is studied comprehensively. Let $C_{\alpha}$ denote the standard 
$\alpha$-dimensional Cantor set contained in the interval $[0,1]$. Let 

$$F_{\alpha}=(C_{\alpha}-1) \cup (C_{\alpha}+1)$$ and let $\mu$ denote the natural measure on this set. Let $E=F_{\alpha} \times F_{\beta}$. Observe that we can a fit a 
$\sqrt{\epsilon}$ by $\epsilon$ rectangle in the annulus $\{x: 1 \leq |x| \leq 1+\epsilon \}$ near $(0, \pm 1)$ and also near $(\pm 1, 0)$. 

Fix $x$ and observe that 

$$ (\mu \times \mu) \{(y,z): 1 \leq |x-z| \leq 1+\epsilon; 1 \leq |x-y| \leq 1+\epsilon; \sqrt{2} \leq |y-z| \leq 
\sqrt{2}+\epsilon \}$$

$$ \approx \epsilon^{\frac{\alpha}{2}+\beta} \cdot \epsilon^{\alpha+\beta}=
\epsilon^{\frac{3}{2} \alpha+2 \beta}.$$ Integrating in $x$, we see that 
$$ (\mu \times \mu \times \mu) \{(x,y,z): 1 \leq |x-z| \leq 1+\epsilon; 1 \leq |x-y| \leq 1+\epsilon; \sqrt{2} \leq |y-z| \leq \sqrt{2}+\epsilon \} \gtrsim \epsilon^{\frac{3}{2} \alpha+2 \beta}.$$ 
We need this quantity to be $\lesssim \epsilon^3$, which leads to the equation 
$$ \frac{3}{2} \alpha+2 \beta \ge 3.$$ 

Choosing $\alpha=1$ and $\beta=\frac{3}{4}$ shows that the inequality (\ref{mamaexpression}) does not in general hold if $s<\frac{7}{4}$. It is important to note that this does not prove that ${\mathcal L}^3(T_2(E))>0$ does not in general hold if $s<\frac{7}{4}$. 

We stress that the calculation above pertains to the trilinear expression (\ref{mamaexpression}). We do not know of any example that shows that ${\mathcal L}^3(T_2(E))$ is not in general positive if the Hausdorff dimension of $E$ is greater than one. The discrepancy here is not particularly surprising because it already takes place in the study of distance sets. For example, as we point out in the introduction, it is known that if the Hausdorff dimension of $E \subset {\Bbb R}^2$ is $\leq 1$, then it is not in general true that ${\mathcal L}^1(\Delta(E))>0$. A result due to Wolff \cite{W99} says that if the Hausdorff dimension of $E$ is greater than $\frac{4}{3}$, then ${\mathcal L}^1(\Delta(E))>0$. On the other hand, an example due to Mattila \cite{Mat87} shows that if the Hausdorff dimension of $E$ is less than $\frac{3}{2}$ and $\mu$ is a Frostman measure on $E$, then  

\begin{equation} \label{falconer} \limsup_{\epsilon \to 0} \epsilon^{-1} (\mu \times \mu) \{(x,y) \in E \times E: 1 \leq |x-y| \leq 1+\epsilon \}=\infty.\end{equation} 

We note that (\ref{falconer}) is the analogue of (\ref{nu}). It says that the distance measure, defined by 
$$ \int f(t) d\nu(t)=\int \int f(|x-y|)\, d\mu(x)\,  d\mu(y),$$ 
{\ag has Radon-Nikodym derivative which is not in $L^{\infty}$.} 

\vskip.125in

\vskip.125in 

\section{Application to discrete geometry} 
\label{application} 

\vskip.125in 

\begin{definition} Let $P$ be a set of $n$ points contained in ${[0,1]}^2$. Define the measure
\begin{equation} \label{pizdatayamera} d \mu^s_P(x)=n^{-1} \cdot n^{\frac{d}{s}} \cdot \sum_{p \in P} \chi_{B^{\, p}_{n^{-\frac{1}{s}}}}(x)\, dx, \end{equation} where $\chi_{B^{\, p}_{n^{-\frac{1}{s}}}}(x)$ is the characteristic function of the ball of radius $n^{-\frac{1}{s}}$ centered at $p$. We say that $P$ is \emph{$s$-adaptable} if 

$$ I_s(\mu_P)=\int \int {|x-y|}^{-s} d\mu^s_P(x) d\mu^s_P(y)<\infty.$$ 

\end{definition} 

\vskip.125in

This is equivalent to the statement 

\begin{equation} \label{discreteenergy} n^{-2} \sum_{p \not=p' \in P} {|p-p'|}^{-s} \lesssim 1.\end{equation} 

\vskip.125in 

To understand this condition in  clearer geometric terms, suppose that $P$ comes from a \linebreak$1$-separated set $A$, scaled down by its diameter. Then the condition (\ref{discreteenergy}) takes the form 

\begin{equation} \label{discreteenergylarge} n^{-2} \sum_{a \not=a' \in A} {|a-a'|}^{-s} \lesssim {(diameter(A))}^{-s}. \end{equation} 
This says $P$ is $s$-adaptable if it is a scaled $1$-separated set where the expected value of the distance between two points raised to the power $-s$ is comparable to the value of the diameter raised to the power of $-s$. This basically means that for the set to be $s$-adaptable, clustering is not allowed to be too severe. 

To put it in more technical terms, $s$-adaptability means that a discrete point set $P$ can be thickened into a set which is uniformly $s$-dimensional in the sense that its energy integral of order $s$ is finite. Unfortunately, it is shown in \cite{IRU10} that there exist finite point sets which are not $s$-adaptable for certain ranges of the parameter $s$. The point is that the notion of Hausdorff dimension is much more subtle than the simple ``size" estimate. However, many natural classes of sets are $s$-adaptable. For example, homogeneous sets studied by Solymosi and Vu \cite{SV04} and others are $s$-adaptable for all $0<s<d$. See also \cite{IJL09} where $s$-adaptability of homogeneous sets is used to extract discrete incidence theorems from Fourier type bounds. 

Before we state the discrete result that follows from Theorem \ref{main}, let us briefly review what is known. If $P$ is set of $n$ points in ${[0,1]}^2$, let $u_{2,2}(n)$ denote the number of times a fixed triangle can arise among points of $P$. It is not hard to see that 
\begin{equation} \label{st} u_{2,2}(n)=O(n^{\frac{4}{3}}).\end{equation} 
This follows easily from the fact that a single distance cannot arise more than $O(n^{\frac{4}{3}})$ times, which, in turn, follows from the celebrated Szemeredi-Trotter incidence theorem. See \cite{BMP05} and the references  therein. By the pigeon-hole principle, one can conclude that 
\begin{equation} \label{distance} \# T_2(P) \gtrsim \frac{n^3}{n^{\frac{4}{3}}}=n^{\frac{5}{3}}.\end{equation}

However, it is not difficult to see that one can do quite a bit better as far as the lower bound on $\# T_2(P)$ is concerned. It is shown in \cite[p. 263]{BMP05} that 
$$ \# T_2(P) \gtrsim n \cdot \# \{|x-y|: x,y \in P\}.$$ 

Guth and Katz have recently settled the Erd\H os distance conjecture in a remarkable paper (\cite{GK10}), proving that 
$$ \# \{|x-y|: x,y \in P\} \gtrsim \frac{n}{\log(n)},$$ and it follows that 

$$ \# T_2(P) \gtrsim \frac{n^2}{\log(n)},$$ which, up to logarithmic factors, is the optimal bound. However, Theorem \ref{main} does allow us to obtain an upper bound on $u_{2,2}$ for $s$-adaptable sets that is better than the one in (\ref{st}). Before we state the main result of this section, we need the following definition. 

\begin{definition} \label{almost} Let $P$ be a subset of ${[0,1]}^2$ consisting of $n$ points as. Let $\delta>0$ and define 
$$u^{\delta}_{2,2}(n)=\# \{(x^1,x^2,x^3) \in P \times P \times P: 
t_{ij}-\delta \leq |x^i-x^j| \leq t_{ij}+\delta \},$$ where the dependence on $t=\{t_{ij}\}$ is suppressed. 
\end{definition} 

Observe that obtaining an upper bound for $u^{\delta}_{2,2}(n)$ with arbitrary $t_{ij}$ immediately
 implies the same upper bound on $u_{2,2}(n)$ defined above. The main result of this section is the following. 

\begin{corollary} \label{discreteresult} Suppose $P \subset {[0,1]}^2$ is $s$-adaptable for $s=\frac{7}{4}+a$ for every sufficiently small $a>0$. Then for every $b>0$, there exists $C_{b}>0$ such that 
\begin{equation} \label{discretekaifest} u^{n^{-\frac{4}{7}-b}}_{2,2}(n) \leq C_{b} 
n^{\frac{9}{7}+b}. \end{equation} 

\end{corollary} 

The proof follows from Theorem \ref{main} in the following way. Let $E$ denote the support of $d\mu^s_P$,  defined as in (\ref{pizdatayamera}) above. We know that if $s>\frac{7}{4}$, then 
\begin{equation} \label{konchayem} (\mu_P^s \times \mu_P^s \times \mu_P^s) \{(x^1,x^2,x^3): t_{ij} \leq |x^i-x^j| \leq t_{ij}+\epsilon \} \lesssim \epsilon^3. \end{equation} 

Taking $\epsilon=n^{-\frac{1}{s}}$, we see that the left hand side is 
$$ \approx n^{-3} \cdot u_{2,2}^{n^{-\frac{1}{s}}}(n)$$ and we conclude that 
$$ u_{2,2}^{n^{-\frac{1}{s}}}(n) \lesssim n^{3-\frac{3}{s}},$$ which yields the desired result since $s=\frac{7}{4}+a$.

As we note above, this result is stronger than the previously known $u_{2,2}(n) \lesssim n^{\frac{4}{3}}$. However, our result holds under an additional restriction that $P$ is $s$-adaptable. We hope to address this issue in a subsequent paper. 
\bigskip


\end{document}